\documentclass[%letterpaper,
twoside, a4paper,
11pt,
]{amsart}

\def\publname{\scriptsize 20/09/2005}
\title[Geometric cohomology frames]{Geometric cohomology frames on Hausmann--Holm--Puppe conjugation spaces}

\copyrightinfo{MMV}{JvH}
\def\currentvolume{}\def\currentissue{}\pagespan{1}{1}\PII{}

\author{Joost van Hamel}
\address{
K.U.Leuven \\
Departement Wiskunde \\
Celestijnenlaan 200B \\
B-3001 Leuven (Heverlee) \\
Belgium 
}
\email{vanhamel@member.ams.org}
\usepackage[shortalphabetic]{amsrefs}

\usepackage{mathptmx}

\theoremstyle{plain}

\newtheorem{theorem}{Theorem}[section]

\newtheorem{lemma}[theorem]{Lemma}
\newtheorem{corollary}[theorem]{Corollary}

\newtheorem*{theorem*}{Theorem}
\newtheorem*{maintheorem}{Theorem}
\newtheorem*{proposition*}{Proposition}
\newtheorem*{lemma*}{Lemma}
\newtheorem*{corollary*}{Corollary}

\theoremstyle{remark}
\newtheorem{remark}[theorem]{Remark}
\newtheorem{Remarks}[theorem]{Remarks}

\newtheorem*{acknowledgements}{Acknowledgements}

\newenvironment{remarks}{\begin{Remarks}\nopagebreak[4]
\rule{1em}{0ex}\par % "\emptyline"
\begin{theoremlist}}%
{\end{theoremlist}\end{Remarks}}

\newcounter{thlist}
\newenvironment{theoremlist}{\begin{list}{\hfill\textup{(\roman{thlist})}}%
{\usecounter{thlist}
\setlength{\leftmargin}{0pt}
\setlength{\rightmargin}{0pt}
\setlength{\labelwidth}{-1.5em}
}}{\end{list}}

\newcounter{assum}
\renewcommand{\theassum}{\Alph{assum}}

\newenvironment{assumptionlist}%
  {\begin{list}
      {\theassum}
      {\usecounter{assum}
       
       \setlength{\leftmargin}{2em}%{0pt}
        \setlength{\rightmargin}{0pt}%{0pt}
        \setlength{\labelwidth}{1.5em}%{-1.5em}}}
        \setlength{\labelsep}{.5em}
       }
  }
  {\end{list}}

\renewcommand{\labelenumi}{(\roman{enumi})}

\newtheorem*{remark*}{Remark}

\newcommand{\mathb}[1]{\mathbf{#1}}
 \newcommand{\fcolon}{\colon\mskip\thinmuskip} %the command giving the colon in 
 \newcommand{\scolon}{\colon\mskip\thinmuskip} %the command giving the colon in 
 \newcommand{\isoto}{\overset{\sim}{\to}}
 \newcommand{\into}{\hookrightarrow}
 
 \DeclareMathOperator{\identity}{id}

 \newcommand{\Rr}{{\mathb{R}}}
 \newcommand{\Cc}{{\mathb{C}}}
 \newcommand{\Ff}{{\mathb{F}}}
 
 \newcommand{\hvar}{{-}}

\newcommand{\sH}{\mathcal{H}}
\newcommand{\sZ}{\mathcal{Z}}

\renewcommand{\phi}{\varphi}

\newcommand{\Ctwo}{{C_2}}

\begin{document}

%%%%%%%%%%%%%%%%%%%%%%%%%%%%%%%%%%%%%%%%%%%%%%%%%%%%%%%%%%%%%%%%%%%%%%%%

\begin{abstract}
For certain manifolds with an involution the
mod 2 cohomology ring of the set of fixed points is isomorphic to the
cohomology ring of the manifold, up to dividing the degrees by two. 
Examples include complex
projective spaces and Grassmannians with the standard antiholomorphic
involution (with real projective spaces and Grassmannians as fixed
point sets).

Hausmann, Holm and Puppe have put this observation in the framework of
equivariant cohomology, and come up with the concept of
\emph{conjugation spaces}, where the ring homomorphisms arise
naturally from the existence of what they call \emph{cohomology frames}. 
Much earlier, Borel and Haefliger had studied the degree-halving
isomorphism between the cohomology rings of complex and real
projective spaces and Grassmannians using the theory of complex and
real analytic cycles and cycle maps into cohomology.

The main result in the present note gives a (purely topological)
connection between these two results and provides a geometric
intuition between the concept of a cohomology frame. 
In particular, we see that if 
every cohomology class on a manifold $X$ with involution
is the Thom class of an
equivariant topological cycle of codimension twice the codimension of
its fixed points (inside the fixed point set of $X$),
these topological cycles will give rise to a cohomology frame.
\end{abstract}

%%%%%%%%%%%%%%%%%%%%%%%%%%%%%%%%%%%%%%%%%%%%%%%%%%%%%%%%%%%%%%%%%%%%%%%%

\maketitle

Let $X$ be a topological space with a continuous involution $\tau$.
We denote the cyclic group of order two by $\Ctwo$.
Since the cohomology ring of $B\Ctwo = P^\infty(\Rr)$, with coefficients in
the field $\Ff_2$ with two elements, is isomorphic to a polynomial ring
$\Ff_2[u]$ in one variable (with $u$ of degree $1$), the restriction
map from $X$ to the fixed point set $X^\tau$ 
in Borel's version of equivariant
cohomology
$H^*_{\Ctwo}(\hvar; \Ff_2):=H^*(E\Ctwo\times_\Ctwo\hvar; \Ff_2)$
(\cite{Borel:trafo})
gives a homomorphism of $\Ff_2[u]$-algebras
\begin{equation*}
r \fcolon H^*_\Ctwo(X; \Ff_2) \to 
H^*_\Ctwo(X^\tau; \Ff_2)[u] = H^*(X^\tau; \Ff_2) \otimes_{\Ff_2} \Ff_2[u].
\end{equation*}
The localisation theorem in equivariant cohomology tells us that if
$X$ is finite-dimensional, this homomorphism becomes 
an isomorphism when we invert the variable $u$.
Together with the ring homomorphism
\begin{equation*}
\rho\fcolon H^*_\Ctwo(X; \Ff_2) \to H^*(X; \Ff_2)
\end{equation*}
this provides a close, but somewhat indirect relation between the
cohomology of $X$ and the cohomology of $X^\tau$.

In the paper~\cite{H-H-P:conjugation}
J.-C.~Hausmann, T.~Holm and V.~Puppe have exhibited an interesting class of
spaces with an involution where this relation can be lifted to a ring
homomorphism between the cohomology of $X$ and the cohomology of
$X^\tau$ which divides the degrees by two.

They call these spaces ``conjugation spaces'', and they are defined to
be those spaces $X$ with an involution $\tau$ for which
$H^\mathrm{odd}(X; \Ff_2) = 0$, and which admit what they call a
\emph{cohomology frame}: a pair $(\kappa, \sigma)$ of additive
homomorphisms
\begin{align*}
\kappa\fcolon H^{2*}(X; \Ff_2) & \to H^*(X^\tau; \Ff_2), \\
%\intertext{and}
\sigma\fcolon H^{2*}(X; \Ff_2) & \to H^{2*}_\Ctwo(X; \Ff_2), 
\end{align*}
such that $\kappa$ is a degree-halving isomorphism,
$\rho \circ \sigma = \identity$, and for every
$m \geq 0$ and every $a \in H^{2m}(X; \Ff_2)$ we have 
the so-called \emph{conjugation relation}
\begin{equation}\label{eq:conj-rel}
r \circ \sigma(a) = \kappa(a) u^m + 
\omega_{m+1}u^{m-1} + \omega_{m+2} u^{m-2} + \cdots +
\omega_{2m}
\end{equation}
with $\omega_i \in H^i(X^\tau; \Ff_2)$ for $i = m+1, \dots 2m$.
What is remarkable about this definition is that the conjugation
relation implies that
\begin{itemize}
\item such a pair $(\kappa, \sigma)$ is unique (if it exists),
\item the homomorphisms $\kappa$ and $\sigma$ are ring homorphisms.
\end{itemize}

The main example in~\cite{H-H-P:conjugation} of a conjugation space is
a so-called \emph{spherical conjugation
  complex}, which is constructed by attaching 
\emph{conjugation cells}: closed
unit disks in $\Cc^n$ with the involution corresponding to complex
conjugation. 
Simple examples of spherical conjugation complexes are complex
projective spaces and complex Grassmannians with the involution given
by complex conjugation (hence the fixed point sets are real projecive
spaces and real Grassmannians).

Much earlier, A.~Borel and A.~Haefliger had studied the degree-halving
isomorphism between the cohomology rings of complex and real
projective spaces and Grassmannians from a different point of view
(and without using equivariant cohomology).
Namely, it follows from Proposition~5.15 in their classic
paper~\cite{B-H} that if for
a compact complex analytic variety $X$ with an antiholomorphic
involution~$\tau$ for which
the $\tau$-invariant analytic cycles 
generate the cohomology $H^*(X; \Ff_2)$ 
and the fixed points of these analytic cycles (which are real
analytic cycles) generate the cohomology
$H^*(X^\tau; \Ff_2)$ of the fixed points, then
the cycle maps induce the desired isomorphism of rings $H^{2*}(X,
\Ff_2) \isoto H^*(X^\tau; \Ff_2)$.

The aim of this note is to elaborate on the connection between these
two points of view, and to provide a geometric insight in the
conjugation equation~\eqref{eq:conj-rel}. 
We will see that for a complex analytic variety as above, 
the $\tau$-equivariant analytic cycles give a geometric construction
of a cohomology frame $(\kappa, \sigma)$, where the conjugation
relation is satisfied because taking fixed points for the
antiholomorphic involution halves the dimension of both the ambient
space and the analytic cycles.

We will work in a more
general topological framework, where we replace complex analytic manifolds
with an antiholomorphic involution
and equivariant analytic cycles by topological analogues.
Since analytic cycles are sums of analytic
subvarieties that can be singular (e.g., Schubert cells), it is
appropriate for a topological generalisation 
to work with singular topological varieties in the sense
of \cite{B-H}. See Section~\ref{sec:thom} for basic definitions and a
discussion of the cohomology classes $[Z] \subset H^{*}(X,
\Ff_2)$, $[Z^\tau] \subset H^*(X^\tau; \Ff_2)$ and $[Z]_\Ctwo \subset
H^*_\Ctwo(X; \Ff_2)$ that may be associated to an equivariant
singular topological subvariety $Z \subset X$ of a topological
manifold $X$ with a locally linear involution.

\begin{maintheorem}
Let $X$ be a (not necessarily compact) connected
topological manifold
with a locally linear involution~$\tau$.
If:
\begin{assumptionlist}
\item
$H^\textrm{odd}(X; \Ff_2) = 0$,
\item
for every $k \geq 0$ 
we have a set $\sZ^k$  of good equivariant singular topological 
subvarieties $Z \subset X$ of codimension~type $(2k, k)$ representing
a basis of the cohomology group $H^{2k}(X; \Ff_2)$,
\item
\label{item:kappa-inj-surj}
either:\\
for every $k \geq 0$ the classes represented by the fixed point
sets of the $Z \in \sZ^k$
are linearly independent in $H^{k}(X^\tau; \Ff_2)$,\\
or:\\ 
for every $k \geq 0$ the classes represented by the fixed point
sets of the $Z \in \sZ^k$ generate $H^{k}(X^\tau; \Ff_2)$,
\end{assumptionlist}
then the homomorphisms
\begin{align*}
\kappa\fcolon H^{2*}(X; \Ff_2) & \to H^*(X^\tau; \Ff_2), &
%\text{and} &&
\sigma\fcolon H^{2*}(X; \Ff_2) & \to H^{2*}_\Ctwo(X; \Ff_2), 
\intertext{defined on the basis elements $\{ [Z] \scolon Z \in \sZ^*\}$ by}
\kappa([Z]) & = [Z^\tau], &
%&&
\sigma([Z]) &= [Z]_\Ctwo,
\end{align*}
form a cohomology frame for $X$.

Moreover, if $X$ is compact, then the first two conditions imply the
third.
\end{maintheorem}

The terminology of singular topological varieties and the cohomologiy
classes they represent is explained in Section~\ref{sec:thom}. This
section also contains the key technical results
Lemma~\ref{lem:restr-thom-class-smooth} and
Corollary~\ref{cor:restr-thom-class}.  The theorem is proved in
Section~\ref{sec:proof-main}.  

Note that for a
compact complex analytic manifold with an anti-holomorphic involution, the
theorem gives in fact not only a topological explanation for
\cite{B-H}*{Prop.~5.15}, but it gives
a slight strengthening as well, 
since we no longer have to check surjectivity of the cycle map onto
$H^*(X^\tau; \Ff_2)$.

\begin{corollary*}
Let $X$ be a compact complex analytic manifold with an anti-holomorphic
involution~$\tau$. If every class in $H^*(X; \Ff_2)$ can be
represented by a
$\tau$-invariant complex analytic cycle, then:
\begin{theoremlist}
\item
 every class in 
$H^*(X^\tau; \Ff_2)$ can be represented by the fixed points of a
$\tau$-invariant complex analytic cycle,
\item
$X$ admits a cohomology frame.
\end{theoremlist}
\end{corollary*}

\begin{remark*}
In the algebro-geometric context (i.e., $X$ projective), a statement
equivalent to the corollary
was proved by V.A.~Krasnov using 
quite different, not purely topological 
methods (\cite{Krasnov:real-alg-max}*{Th.~0.1}).
\end{remark*}

\begin{acknowledgements}
This work has it origins in discussions I had with Jean-Claude
Hausmann and Tara Holm at MSRI, Berkeley in February 2004 about their
work in progress, and it has benefited from later discussions with
Jean-Claude Hausmann during a visit to Geneva.  I would like to thank
them both for the inspiring discussions, and I would like to
acknowledge the financial support received from MSRI from and the
Swiss National Science Foundation for my visits to Berkeley and
Geneva.

I would also like to thank Volker Puppe for his useful comments on a
preliminary version of this note.
\end{acknowledgements}

\section{Thom classes and singular subvarieties}\label{sec:thom}

Let $X$ be a (not necessarily compact) topological
manifold. In this paper, we adhere to the convention that
manifolds are finite-dimensional and do not have boundaries. 
Let $Z \subset X$ be a closed subspace. We write $H^*(X, X\setminus Z; \Ff_2)$
for the cohomology relative to the complement of $Z$. This cohomology
is also known as cohomology with supports in $Z$.
Recall, that if for each connected component $X_i \subset X$ of
dimension $d_i$ the subspace 
$Z_i:=Z \cap X_i$ has cohomological dimension~$\leq d_i-k$, then $H^m(X,
X\setminus Z; \Ff_2) = 0$ for $m < k$ by Alexander--Lefschetz
duality.
 
In analogy with \cite{B-H}*{2.1} we
say that a closed subspace $Z \subset X$ is a \emph{singular subvariety} of
codimension~$k$ if $Z$ contains a nonempty open topological submanifold $U
\subset Z$ such that the submanifold $U_i :=U \cap X_i$ has codimension $k$ in
each connected component $X_i \subset X$ where $U_i$ is
nonempty, 
and such that each 
$\Sigma_i := Z_i \setminus U_i$ has cohomological dimension
strictly less than $d_i - k$. We call any such $U$ a \emph{fat nonsingular
open} of~$Z$. Note that (as in \cite{B-H}) the open $U$ is not required
to be dense in~$Z$ (and for this reason we are using a liberal, rather
than literal translation  of the original
French term `\'epais'). 
Also note that we do not require $U$ to be a
maximal submanifold, so the corresponding `singular subset'  
$\Sigma := Z \setminus U$ may include points where $Z$ 
is locally a manifold.

Recall, that the \emph{Thom class} $[U] \in H^*(Y; \Ff_2)$ of a closed
submanifold of codimension $k$ in a manifold $Y$ can be characterised
by the fact that it is the unique homogeneous class in $H^k(Y, Y\setminus U;
\Ff_2)$ such that cup-product with $[U]$ induces an isomorphism
\begin{equation*}%\label{eq:thom-iso}
[U] \cup \hvar \fcolon H^k(U; \Ff_2) \isoto H^k(Y, Y\setminus U; \Ff_2).
\end{equation*}

We say that a singular subvariety $Z \subset X$ of codimension $k$
\emph{admits a Thom class} if there is a cohomology class
\begin{equation*}
[Z] \in H^k(X, X\setminus Z; \Ff_2)
\end{equation*}
that maps to the Thom class $[U]$ 
of a fat nonsingular open 
$U$ under the restriction map
\begin{equation*}
H^k(X, X\setminus Z; \Ff_2) \to H^k(X\setminus \Sigma, X\setminus Z; \Ff_2).
\end{equation*}
As in \cite{B-H}*{Prop.~2.3} we see that $[Z]$ is unique if it exists
and that the existence of $[Z]$ is guaranteed if $Z$ admits a fat
nonsingular open $U$ such that each $\Sigma_i \subset X_i$ has
cohomological dimension~$< d_i-k-1$. An example of a singular subvariety
without a Thom class is a bounded closed line segment in the plane. A typical
unbounded example would be the union of three half-lines meeting in
one point.

By abuse of notation, we will also write $[Z]$ for the image of the
Thom class in $H^k(X; \Ff_2)$ under the natural map
$H^k(X, X\setminus Z; \Ff_2) \to H^k(X; \Ff_2)$.
If we say (for example, in the statement of the Main Theorem) 
that \emph{$Z$ represents a cohomology class $a \in H^k(X;
\Ff_2)$}, then this combines two assertions: first, that $Z \subset X$ admits a
Thom class, and second, that $a = [Z]$.

\begin{remark}
  For sake of concreteness I have chosen to state and discuss
  everything in terms of \emph{topological} manifolds, singular
  topological subvarieties and locally linear involutions.  However,
  all results and all proofs in this note are valid for
  \emph{cohomological} $\Ff_2$-manifolds, arbitrary continuous
  involutions and the corresponding class of singular subvarieties as
  well.
\end{remark}

\subsection{Equivariant singular subvarieties and equivariant Thom classes}

Assume that the ambient manifold
$X$ admits an involution~$\tau$.
Since for concreteness we work with topological manifolds, we will
assume this involution to be locally linear, so that $X^\tau \subset
X$ is a submanifold.
For a singular topological subvariety $Z\subset X$ of codimension $k$
to be an \emph{equivariant singular topological subvariety}, 
we  require not only that $\tau(Z) = Z$, but also that $Z$ admits a
fat nonsingular open $U \subset Z$ 
such that $\tau$ acts locally linear on the pair
$U \subset X$. This
ensures that $U^\tau \subset X^\tau$ is a submanifold, 
so we call $U$ a \emph{good equivariant submanifold}.
Of
course, this condition on $U$ 
does not ensure that $Z^\tau$ is a singular subvariety with
fat nonsingular open $U^\tau$, since the dimension of $\Sigma^\tau$
may be too big. We say that $Z \subset X$ is a
\emph{good equivariant singular topological subvariety of codimension
  type $(k, k')$} if $Z$ admits a fat nonsingular open $U \subset Z$
  such that $Z^\tau$
  is a singular
  topological subvariety of codimension $k'$ in $X^\tau$ with fat
  nonsingular open~$U^\tau$.

Writing\begin{equation*}
H^*_\Ctwo(X, X\setminus Z; \Ff_2) := H^*(E\Ctwo \times_\Ctwo X,
E\Ctwo \times_\Ctwo (X\setminus Z);\Ff_2)
\end{equation*}
for the equivariant cohomology with supports in $Z$, the Borel--Serre
spectral sequence
\begin{equation*}
E_2^{pq} = H^p(B\Ctwo, \sH^q(X, X\setminus Z; \Ff_2)) \Rightarrow 
H^{p+q}_\Ctwo(X,X\setminus Z; \Ff_2),
\end{equation*}
(i.e., the Leray spectral sequence for the fibration 
$E\Ctwo \times_\Ctwo X \to B \Ctwo$)
tells us that $H^j_\Ctwo(X, X\setminus Z; \Ff_2) = 0$ when $d-j$ is
greater than the
cohomological dimension of $Z$ and that
\begin{equation}\label{eq:inj-rho}
\rho \fcolon H^j_\Ctwo(X, X\setminus Z; \Ff_2) \to H^j(X, X\setminus Z; \Ff_2)
\end{equation}
is injective when $d-j$ is equal to the cohomological dimension of $Z$, with
the image equal to the $\Ctwo$-invariant cohomology classes.
It follows that if $Z$ admits a nonequivariant Thom class 
$[Z] \in H^k(X, X\setminus Z; \Ff_2)$, then $[Z]$ lifts to a unique class
\begin{equation*}
[Z]_\Ctwo \in H^k_\Ctwo(X, X\setminus Z; \Ff_2)
\end{equation*}
which we call the \emph{equivariant Thom class}.

The equivariant Thom class 
$[U]_\Ctwo \in H^*_\Ctwo(Y, Y\setminus U; \Ff_2)$ 
of a closed equivariant submanifold $U$ of a manifold with involution
$Y$ can be
characterised by the fact that it is the unique homogeneous class in
$H^*_\Ctwo(Y, Y\setminus U; \Ff_2)$ such that cup-product with $[U]_\Ctwo$
induces an isomorphism
\begin{equation}\label{eq:equiv-thom-iso}
[U]_\Ctwo \cup \hvar \fcolon H^*_\Ctwo(U; \Ff_2) \isoto 
H^*_\Ctwo(Y, Y\setminus U; \Ff_2).
\end{equation}
The equivariant Thom class $[Z]_\Ctwo \in H^*_\Ctwo(X, X\setminus Z;
\Ff_2)$
of a good equivariant
singular topological subvariety $Z \subset X$ can be characterised
(if it exists) as the unique homogeneous class in 
$H^*_\Ctwo(X, X\setminus Z;
\Ff_2)$ that 
maps by restriction to the equivariant Thom class of any fat
equivariant nonsingular open $U$ of $Z$.

The following results links the equivariant Thom class to the Thom class
of the set of fixed points.

\begin{lemma}\label{lem:restr-thom-class-smooth}
Let $U \subset Y$ be a good equivariant closed submanifold of
codimension~$k$ of a not necessarily compact
topological manifold $Y$ with a locally linear involution $\tau$.
Assume that $U^\tau \subset Y^\tau$ has
codimension~$k'$. 
Then the restriction map
\begin{equation*}
r \fcolon H^*_\Ctwo(Y, Y\setminus U; \Ff_2) \to 
H^*_\Ctwo(Y^\tau,
Y^\tau\setminus U^\tau; \Ff_2) = 
H^i(Y^\tau, Y^\tau\setminus U^\tau; \Ff_2)[u]
\end{equation*}
maps the equivariant Thom class 
$[U]_\Ctwo \in H^*_\Ctwo(Y, Y\setminus U; \Ff_2)$ to the class
\begin{equation*}
  r([U]_\Ctwo) = [U^\tau] u^{k-k'} + \eta_{k'+1} u^{k-k'-1} + \dots + 
  \eta_{k}
\end{equation*}
with $\eta_i \in H^i(Y^\tau, Y^\tau\setminus U^\tau; \Ff_2)$, and
$[U^\tau] \in H^{k'}(Y^\tau, Y^\tau\setminus U^\tau; \Ff_2)$ the Thom
class of $U^\tau \subset X^\tau$.
\end{lemma}
\begin{proof}
Since $H^i(Y^\tau, Y^\tau\setminus U^\tau; \Ff_2) = 0$ for $i < k'$ we have
that
\begin{equation*}
  r([U]_\Ctwo) = \eta_{k'} u^{k-k'} + \eta_{k'+1} u^{k-k'-1} + \dots + 
  \eta_{k}.
\end{equation*}
The isomorphism~\eqref{eq:equiv-thom-iso} together with the
localisation theorem implies that cup-product with the class
$r([Z]_\Ctwo) \in H^k_\Ctwo(Y^\tau, Y^\tau\setminus U^\tau; \Ff_2)$
induces an isomorphism
\begin{equation}\label{eq:equiv-thom-iso-r}
  r([U]_\Ctwo) \cup \hvar \fcolon H^*(U^\tau; \Ff_2)[u,u^{-1}] \isoto H^*(Y^\tau, Y^\tau\setminus U^\tau; \Ff_2)[u,u^{-1}].
\end{equation}
It follows that cup-product with $\eta_{k'} \in H^{k'}(Y^\tau,
Y^\tau\setminus U^\tau; \Ff_2)$ induces an isomorphism
\begin{equation*}
  \eta_{k'} \cup \hvar \fcolon H^*(U^\tau; \Ff_2) \isoto 
      H^{*+k'}(Y^\tau, Y^\tau\setminus U^\tau; \Ff_2),
\end{equation*}
hence $\eta_{k'}$ is the Thom class of $U^\tau$.
\end{proof}

\begin{corollary}\label{cor:restr-thom-class}
  Let $Z \subset X$ be a good equivariant singular topological subvariety
  of codimension type $(k, k')$ of a not necessarily compact topological manifold
  $X$ with a locally linear involution $\tau$.

 If $Z \subset X$ admits a Thom class $[Z] \in H^k(X, X\setminus Z;
\Ff_2)$, then 
\begin{theoremlist}
\item
$Z \subset X$ admits an
equivariant Thom class $[Z]_\Ctwo \in H^k_\Ctwo(X, X\setminus Z;
\Ff_2)$
and the natural map
\begin{equation*}
\rho\fcolon H^k_\Ctwo(X, X\setminus Z; \Ff_2) \to H^k(X, X \setminus Z; \Ff_2)
\end{equation*}
sends $[Z]_\Ctwo$ to $[Z]$,
\item
$Z^\tau \subset X^\tau$ admits a Thom class $[Z^\tau] \in
H^{k'}(X^\tau, X \setminus Z^\tau; \Ff-2)$, and
the restriction map
\begin{equation*}
r \fcolon H^*_\Ctwo(X, X\setminus Z; \Ff_2) \to 
H^*_\Ctwo(X^\tau,
X^\tau\setminus Z^\tau; \Ff_2) = 
H^*(X^\tau, X^\tau\setminus Z^\tau; \Ff_2)[u]
\end{equation*}
sends the equivariant Thom class 
$[Z]_\Ctwo \in H^*_\Ctwo(X, X\setminus Z; \Ff_2)$ to the class
\begin{equation*}
  r([Z]_\Ctwo) = [Z^\tau] u^{k-k'} + \omega_{k'+1} u^{k-k'-1} + \dots + 
  \omega_{k}
\end{equation*}
with $\omega_i \in H^i(X^\tau, X^\tau\setminus Z^\tau; \Ff_2)$, and
$[Z^\tau] \in H^{k'}(X^\tau, X^\tau\setminus Z^\tau; \Ff_2)$ the Thom
class of $Z^\tau \subset X^\tau$.
\end{theoremlist}
\end{corollary}
\begin{proof}
The existence of a unique lift $[Z]_\Ctwo$ of $[Z]$ was already proved
above, so we only have to prove the second part.
Since $Z^\tau \subset X^\tau$ is a singular topological subvariety of
dimension~$k'$, we have that
\begin{equation*}
  r([Z]_\Ctwo) = \omega_{k'} u^{k-k'} + \omega_{k'+1} u^{k-k'-1} + \dots + 
  \omega_{k}
\end{equation*}
with $\omega_i \in H^i(X^\tau, X^\tau\setminus Z^\tau; \Ff_2)$.
Let $j \fcolon U \into Z$ be a nice fat equivariant nonsingular open.
By definition, $j^*([Z]_\Ctwo) = [U]_\Ctwo$. Since $r \circ j^* = j^*
\circ r$, Lemma~\ref{lem:restr-thom-class-smooth} implies that
\begin{equation*}
  j^*\omega_{k'} u^{k-k'} + j^*\omega_{k'+1} u^{k-k'-1} + \dots + 
  j^* \omega_{k} = [U^\tau] u^{k-k'} + \eta_{k'+1} u^{k-k'-1} + \dots + 
  \eta_{k}.
\end{equation*}
Hence $\omega_{k'}$ is the Thom class of $Z^\tau \subset X^\tau$.
\end{proof}

\begin{remarks}\label{rems}

\item
A treatment more in the spirit of~\cite{B-H} would 
be in terms of equivariant
Borel--Moore homology and fundamental classes (compare
\cite{vH:CWITract}*{III.7.1, IV.1}), which is linked to the present
treatment via equivariant Poincar\'e duality (\emph{loc.\ cit.}).
The above treatment was chosen to remain closer to the language and
approach in~\cite{H-H-P:conjugation}. 

\item
Observe that \emph{a priori} it is not obvious at all that the existence of a
Thom class for~$Z$ implies the existence of a Thom class
for~$Z^\tau$. 
For example, in~\cite{B-H}, complexifications of real analytic sets
play an important role, but the existence of a fundamental
class for a real analytic set is proved using the highly
non-topological operation of \emph{normalisation}, whereas
for a complex analytic space the existence of a fundamental class 
is deduced from the
simple topological observation that the set of singular points is of
topological codimension~$\geq 2$.

\item The proof of Lemma~\ref{lem:restr-thom-class-smooth} is
  analogous to the proof of
  \cite{vH:CWITract}*{Th.~III.7.4}. This proof was
  inspired by the proof of~\cite{Allday-Puppe}*{Prop.~5.3.7}, which
  proves the case $U = Y^\tau$ of Lemma~\ref{lem:restr-thom-class-smooth} in
  the case where $Y$ is a Poincar\'e duality space, rather than a
  manifold.
\end{remarks}

\section{Proof of the theorem}\label{sec:proof-main}

Let $X$ be a (not necessarily compact) connected
topological manifold
with a locally linear involution~$\tau$.
Assume that
\begin{assumptionlist}
\item
$H^\textrm{odd}(X; \Ff_2) = 0$,
\item
for every $k \geq 0$ 
we have a set $\sZ^k$  of good equivariant singular topological 
subvarieties $Z \subset X$ of codimension type~$(2k,k)$ representing
a basis of the cohomology group $H^{2k}(X; \Ff_2)$,
\end{assumptionlist}

We define homomorphisms
\begin{align*}
\kappa\fcolon H^{2*}(X; \Ff_2) & \to H^*(X^\tau; \Ff_2), &
%\text{and} &&
\sigma\fcolon H^{2*}(X; \Ff_2) & \to H^{2*}_\Ctwo(X; \Ff_2), 
\intertext{on the basis elements $\{ [Z] \scolon Z \in \sZ^*\}$ by}
\kappa([Z]) & = [Z^\tau], &
%&&
\sigma([Z]) &= [Z]_\Ctwo.
\end{align*}
where the classes $[Z]_\Ctwo \in H^*_\Ctwo(X; \Ff_2)$ and 
$[Z^\tau] \in H^*(X^\tau; \Ff_2)$ exist by
Corollary~\ref{cor:restr-thom-class}.

By construction, $\kappa$ is a degree-halving homomorphism, $\rho \circ
\sigma = \identity$, and the conjugation relation~\eqref{eq:conj-rel}
holds by Corollary~\ref{cor:restr-thom-class}.
In order to prove that $\kappa$ is an isomorphism, 
we observe that the existence of $\sigma$ implies that 
$X$ is \emph{equivariantly formal}, i.e., we have an isomorphism of
$\Ff_2[u]$-modules
\begin{equation*}
H^*(X; \Ff_2)[u] \isoto H^*_\Ctwo(X; \Ff_2).
\end{equation*}
It follows from the localisation theorem that
this isomorphism induces an isomorphism
\begin{equation*}
H^*(X; \Ff_2)[u,u^{-1}] \isoto H^*(X^\tau; \Ff_2)[u,u^{-1}].
\end{equation*}
Putting $u = 1$, we get an isomorphism of (possibly infinite
dimensional) $\Ff_2$-vector spaces
\begin{equation}\label{eq:nongrad-iso}
H^*(X; \Ff_2) \isoto H^*(X^\tau; \Ff_2),
\end{equation}
which does \emph{not} preserve (or halve) the grading, but which, 
by the above hypotheses
and~Corollary~\ref{cor:restr-thom-class}, does map $H^{\geq 2k}(X; \Ff_2)$
to $H^{\geq k}(X^\tau; \Ff_2)$ for every $k \geq 0$.
In other words, with the appropriate filtrations on source and
target, 
the isomorphism~\eqref{eq:nongrad-iso} is an
isomorphism of filtered vector spaces that preserves the 
filtrations, with $\kappa$ the corresponding homomorphism of the graded
quotients associated to the filtrations.
Since the filtrations are of finite length,
 $\kappa$ is an isomorphism if it is either injective or
surjective, which is the case by the fourth hypothesis:
\begin{assumptionlist}
\item[\ref{item:kappa-inj-surj}]
either:\\
for every $k \geq 0$ the classes represented by the fixed point
sets of the $Z \in \sZ^k$
are linearly independent in $H^{k}(X^\tau; \Ff_2)$,\\
or:\\ 
for every $k \geq 0$ the classes represented by the fixed point
sets of the $Z \in \sZ^k$ generate $H^{k}(X^\tau; \Ff_2)$
\end{assumptionlist}
This finishes the proof of the first part of the theorem.

For the final assertion of the theorem we let  $X$ be compact, satisfying 
the first two hypotheses.
These hypotheses imply that $X$ is even-dimensional and that
 the dimension~$d$ of $X^\tau$ is half the dimension of $X$. 
We will establish the injectivity of $\kappa$ by
a weak form of equivariant intersection theory.
Assume that we have a $k \geq 0$ and an $a \in H^{2k}(X; \Ff_2)$
such that
$\kappa(a) = 0$. The conjugation equation implies 
that
\begin{equation*}
r \circ \sigma(a) = \omega_{k+1}u^{k-1} + \omega_{k+2} u^{k-2} + \cdots +
\omega_{2k}.
\end{equation*}
with $\omega_i \in H^i(X^\tau; \Ff_2)$.
By Poincar\'e duality we have a cohomology class 
$b \in H^{2d-2k}(X; \Ff_2)$ with
$a \cup b \neq 0$,
hence
\begin{equation*}
\sigma(a) \cup \sigma(b) \neq 0.
\end{equation*} 
On the other hand, the conjugation equation implies
\begin{equation*}
r \circ \sigma(b) = \eta_{d-k} u^{d-k} + \eta_{d-k+1} u ^{d-k-1} +
\cdots + \eta_{2d-2k}
\end{equation*}
with $\eta_i \in H^i(X^\tau; \Ff_2)$ (in fact, $\eta_{d-k} = \kappa(b)$,
but we will not need that here).
Since $\omega_i \cup \eta_j = 0$ for $i+j > d$, we see that
\[r \circ \sigma(a) \cup r \circ \sigma(b) = 0,\]
which implies that $r \circ \sigma (a \cup b) = 0$.
Since $\sigma(a \cup b) = \sigma(a) \cup \sigma(b) \neq 0$, this 
contradicts the injectivity of $r$ 
(which follows from the localisation theorem).
Hence we have shown that $\kappa$ is injective, which finishes the
proof of the theorem. 
\qed

\providecommand{\refjournaltitle}[1]{\textit{#1}}

\end{document}